\documentclass[12pt]{amsart} 
\usepackage{amssymb,amsmath,amscd} 
\theoremstyle{plain} 
\newtheorem{Lem}{Lemma}[section] 
\newtheorem{Prop}[Lem]{Proposition} 
\newtheorem{Thm}[Lem]{Theorem} 
\newtheorem{Cor}[Lem]{Corollary} 
\theoremstyle{definition} 
\newtheorem{Def}[Lem]{Definition} 
\newtheorem{Rem}[Lem]{Remark} 
 
\newtheorem{Pro}[Lem]{Problem} 
\newtheorem{Set}[Lem]{Setting} 
\errorcontextlines=0  \numberwithin{equation}{section}

\newcommand{\bbC}{{\mathbb C}} 
\newcommand{\bbR}{{\mathbb R}}

\newcommand{\bbK}{{\mathbb K}}

\newcommand{\calH}{{\mathcal H}} 
\newcommand{\calL}{{\mathcal L}}

\newcommand{\calS}{{\mathcal S}} 
\newcommand{\calT}{{\mathcal T}} 
\newcommand{\calV}{{\mathcal V}} 
\newcommand{\id}{\mathrm{id}} 
 
\newcommand{\sign}{\mathrm{sign}} 
\newcommand{\frakl}{{\mathfrak l}} 
\newcommand{\frakr}{{\mathfrak r}} 
\newcommand{\fraka}{{\mathfrak a}} 

\newcommand{\fR}{{\mathfrak R}}
 
\newcommand{\bpi}{\begin{picture}} 
\newcommand{\epi}{\end{picture}}

\newcommand{\regrep}{\breve\omega}
\newcommand{\twoline}[2]{\genfrac{}{}{0pt}{}{#1}{#2}}

\begin{document} 

\title[Symmetry classes of covariant derivatives] 
{On the symmetry classes of the first covariant derivatives of tensor fields} 
\date{September 2002} 
\author[B. Fiedler]{Bernd Fiedler}  
\address{Bernd Fiedler \\ Mathematisches Institut \\ Universit\"at Leipzig\\ 
Augustusplatz 10/11 \\ D-04109 Leipzig \\ Germany}
\urladdr{http://home.t-online.de/home/Bernd.Fiedler.RoschStr.Leipzig/}  
\email{Bernd.Fiedler.RoschStr.Leipzig@t-online.de}  
\subjclass{53B20, 15A72, 05E10, 16D60, 05-04} 

\begin{abstract}
We show that the symmetry classes of torsion-free covariant derivatives $\nabla T$ of $r$-times covariant tensor fields $T$ can be characterized by Littlewood-Richardson products $\sigma [1]$ where $\sigma$ is a representation of the symmetric group $\calS_r$ which is connected with the symmetry class of $T$. If $\sigma \sim [\lambda]$ is irreducible then $\sigma [1]$ has a multiplicity free reduction $[\lambda] [1] \sim \sum_{\lambda \subset \mu} [\mu]$ and all primitive idempotents belonging to that sum can be calculated from a generating idempotent $e$ of the symmetry class of $T$ by means of the irreducible characters or of a discrete Fourier transform of $\calS_{r+1}$. We apply these facts to derivatives $\nabla S$, $\nabla A$ of symmetric or alternating tensor fields.
The symmetry classes of the differences $\nabla S - \mathrm{sym}(\nabla S)$ and $\nabla A - \mathrm{alt}(\nabla A) = \nabla A - dA$ are characterized by Young frames
$(r , 1) \vdash r+1$ and $(2 , 1^{r-1}) \vdash r+1$, respectively. However, while the symmetry class of $\nabla A - \mathrm{alt}(\nabla A)$ can be generated by Young symmetrizers of $(2 , 1^{r-1})$, no Young symmetrizer of 
$(r , 1)$ generates the symmetry class of $\nabla S - \mathrm{sym}(\nabla S)$.
Furthermore we show in the case $r = 2$ that $\nabla S - \mathrm{sym}(\nabla S)$ and $\nabla A - \mathrm{alt}(\nabla A)$ can be applied in generator formulas of algebraic covariant derivative curvature tensors.
For certain symbolic calculations we used the Mathematica packages {\sf Ricci} and {\sf PERMS}.
\end{abstract}

\maketitle 

%
%

\section{Introduction}
The investigations of the present paper arose from the search for generators of {\it algebraic curvature tensors}. Algebraic curvature tensors are covariant tensors of order 4 which have the same algebraic properties as the {\it Riemannian curvature tensor}.
\begin{Def}
A covariant tensor $\fR$ of order 4 is called an {\it algebraic curvature tensor} iff its coordinates satisfy the conditions
\begin{eqnarray}
{\fR}_{i j k l} \;=\; - {\fR}_{j i k l} \;=\; - {\fR}_{i j l k} & = & {\fR}_{k l i j} \label{tausch}\\
{\fR}_{i j k l} \,+\, {\fR}_{i k l j} \,+\, {\fR}_{i l j k} & = & 0 \label{bia1}\,.
\end{eqnarray}
A covariant tensor $\fR'$ of order 5 is called an {\it algebraic covariant derivative curvature tensor} iff its coordinates fulfil
\begin{eqnarray}
{\fR}_{i j k l m}' \;=\; - {\fR}_{j i k l m}' \;=\; - {\fR}_{i j l k m}' & = & {\fR}_{k l i j m}' \\
{\fR}_{i j k l m}' \,+\, {\fR}_{i k l j m}' \,+\, {\fR}_{i l j k m}' & = & 0 \\
{\fR}_{i j k l m}' \,+\, {\fR}_{i j l m k}' \,+\, {\fR}_{i j m k l}' & = & 0 \,. \label{bia2}
\end{eqnarray}
\end{Def}
Relation (\ref{tausch}) represents the index commutation symmetry of the Riemannian curvature tensor $R$ whereas relations (\ref{bia1}) and (\ref{bia2}) correspond to the first and second Bianchi identity for the Riemann tensor
\begin{eqnarray*}
R_{i j k l} \,+\, R_{i k l j} \,+\, R_{i l j k} & = & 0\\
R_{i j k l \,;\, m} \,+\, R_{i j l m \,;\, k} \,+\, R_{i j m k \,;\, l} & = & 0 \,.
\end{eqnarray*}
Investigations of algebraic curvature tensors were carried out by many authors. (See the extensive bibliography in the book \cite{gilkey5} by P. B. Gilkey.) One of the problems which are considered in connection with algebraic curvature tensors is the {\it search for generators} of algebraic curvature tensors. P. B. Gilkey \cite[pp.41-44]{gilkey5} and B. Fiedler \cite{fie20} gave different proofs that the vector space of all algebraic curvature tensors is spanned by each of the following types of tensors
\begin{eqnarray}
y_t^{\ast}(S \otimes S) & \;\;\;,\;\;\; & y_t^{\ast}(A \otimes A) \label{rgen}
\end{eqnarray}
which are defined by symmetric or alternating covariant tensors $S$ or $A$ of order 2. The vector space of algebraic covariant derivative curvature tensors is generated by each of the following tensor types
\begin{eqnarray}
y_{t'}^{\ast}(S \otimes S') & \;\;\;,\;\;\; & y_{t'}^{\ast}(S' \otimes S) \nonumber \\
y_{t'}^{\ast}(U \otimes S) & \;\;\;,\;\;\; & y_{t'}^{\ast}(S \otimes U) \label{divgen} \\
y_{t'}^{\ast}(U \otimes A) & \;\;\;,\;\;\; & y_{t'}^{\ast}(A \otimes U) \nonumber
\end{eqnarray}
(see\footnote{A first proof that $y_{t'}^{\ast}(S \otimes S')$ and $y_{t'}^{\ast}(S' \otimes S)$ are generators for $\fR'$ was given by P. B. Gilkey \cite[p.236]{gilkey5}.} B. Fiedler \cite{fie22}). Here $S$, $A$ are again symmetric or alternating tensors of order 2, $S'$ is a symmetric tensor of order 3 and $U$ is a covariant tensor of order 3 from an irreducible symmetry class that belongs to the partition $(2 , 1) \vdash 3$ and is defined by a minimal right ideal different from the right ideal $f \cdot \bbK [\calS_3]$ with the generating idempotent
\begin{eqnarray}
f \;:=\; \left\{ \frac{1}{2}\,(\id - (1 \,3)) - \frac{1}{6}\,y \right\}
& \;\;\;,\;\;\; &
y \;:=\; \sum_{p \in \calS_3} \mathrm{sign}(p)\,p \,. \label{ipf}
\end{eqnarray}
In (\ref{rgen}) and (\ref{divgen}) $y_t$ and $y_{t'}$ denote the Young symmetrizers of the Young tableaux
\begin{eqnarray}
t \;=\;
\begin{array}{|c|c|}
\hline
1 & 3 \\
\hline
2 & 4 \\
\hline
\end{array}
&\;\;\;,\;\;\; &
t' \;=\;
\begin{array}{|c|c|c|c}
\cline{1-3}
1 & 3 & 5 & \\
\cline{1-3}
2 & 4 &\multicolumn{2}{c}{\;\;\;} \\
\cline{1-2}
\end{array}
\,.
\end{eqnarray}
%
%

The generators (\ref{divgen}) lead to the question whether there exist typical tensorial quantities of differential geometry which possess a symmetry of the same type as the above tensors $U$. It can be shown (see Section \ref{sec3}) that the differences
\begin{eqnarray}
\nabla S \;-\; \mathrm{sym}(\nabla S) & \;\;\;,\;\;\; &
\nabla A \;-\; \mathrm{alt}(\nabla A) \;=\; \nabla A \;-\; \mathrm{d}A \label{das}
\end{eqnarray}
between the covariant derivatives $\nabla S$, $\nabla A$ of symmetric/alternating covariant tensor fields $S$, $A$ of order 2 and their symmetrized/anti-symmetrized covariant derivatives $\mathrm{sym}(\nabla S)$, $\mathrm{alt}(\nabla A)$ have such a symmetry\footnote{Note that $\mathrm{alt}(\nabla A)$ is equal to the exterior derivative $\mathrm{d}A$ of the alternating tensor field $A$.}. Furthermore computer calculations by means of the {\sf Mathematica} packages {\sf Ricci} \cite{ricci3} and {\sf PERMS} \cite{fie10} showed that the symmetry class of $\nabla A \;-\; \mathrm{alt}(\nabla A)$ is generated by certain Young symmetrizers, for instance by the Young symmetrizer of the standard tableau
\begin{eqnarray*}
& &
\begin{array}{|c|c|c}
\cline{1-2}
1 & 3 & \\
\cline{1-2}
2 & \multicolumn{2}{c}{\;\;\;} \\
\cline{1-1}
\end{array}
\,.
\end{eqnarray*}
However, the same computation yields the surprising result that no Young symmetrizer of a Young frame $(2\,1) \vdash 3$ generates the symmetry class of $\nabla S \;-\; \mathrm{sym}(\nabla S)$. One goal of the present paper is to find out whether the covariant derivatives of symmetric or alternating tensor fields of order $r > 2$ have such a behaviour, too.

In Section \ref{sec2} we collect some basic facts about the connection between symmetry classes of covariant tensors of order $r$ and left or right ideals of the group ring $\bbK [\calS_r]$ of the symmetric group $\calS_r$.
In Section \ref{sec3} we show that the symmetry class of a torsion-free covariant derivative $\nabla T$ of a differentiable tensor field $T$ of order $r$ is defined by a left ideal $\frakl \subseteq \bbK [\calS_{r+1}]$ which is the representation space of a Littlewood-Richardson product $\sigma [1]$, where $\sigma$ is a representation of $\calS_r$ connected with the symmetry class of $T$. If $\sigma \sim [\lambda]$, $\lambda \vdash r$, is irreducible then the Littlewood-Richardson rule yields a multiplicity-free decomposition
\begin{eqnarray}
[\lambda][1] & \sim & \sum_{\twoline{\mu \vdash r+1}{\lambda \subset \mu}} \,[\mu]\,. \label{lrsum2}
\end{eqnarray}
In particular the symmetry classes of the covariant derivatives $\nabla S$, $\nabla A$ of symmetric/alternating tensor fields $S , A$ of order $r$ are characterized by Littlewood-Richardson products
\begin{eqnarray*}
[r][1] \;\sim\; [r+1] + [r , 1]
& \;\;\;,\;\;\; &
[1^r][1] \;\sim\; [1^{r+1}] + [2 , 1^{r-1}]\,. 
\end{eqnarray*}
If we know a primitive generating idempotent $e \in \bbK [\calS_r]$ of the symmetry class of $T$, then all unique primitive idempotents $h_{\mu} \in \bbK [\calS_{r+1}]$ corresponding to (\ref{lrsum2}) can be calculated from $e$ by means of the symmetrizers of the irreducible characters of $\calS_{r+1}$ or, more efficiently, by a discrete Fourier transform.

In Section \ref{sec4} we investigate the parts $[r , 1]$, $[2 , 1^{r-1}]$ of $\nabla S$, $\nabla A$ for arbitrary order $r \ge 2$. We show for $\nabla S$ that no Young symmetrizer with a Young frame $(r , 1) \vdash r + 1$ is a generator of the $[r , 1]$-part of $\nabla S$. The $[2 , 1^{r-1}]$-part of $\nabla A$, however, is generated by the Young symmetrizer of the lexicographically greatest standard tableau of $(2 , 1^{r-1}) \vdash r + 1$ (see (\ref{reptabl})) and every other standard tableau of $(2 , 1^{r-1}) \vdash r + 1$ annihilates the $[2 , 1^{r-1}]$-part. Furthermore that $[2 , 1^{r-1}]$-part is generated or annihilated by many other Young symmetrizers of non-standard tableaux of $(2 , 1^{r-1}) \vdash r + 1$. We present complete computer generated lists of such Young symmetrizers for $r = 2, 3, 4$ in an Appendix.

The last Section of the paper deals with the question whether tensors (\ref{das}) can be used as generators $U$ of algebraic covariant derivative curvature tensors in formulas (\ref{divgen}). Both $S$ and $A$ satisfy the condition that the symmetry classes of the tensors (\ref{das}) are not generated by the above right ideal $f \cdot \bbK[\calS_3]$ with generating idempotent (\ref{ipf}). Thus all tensors (\ref{das}) can play the role of $U$ in (\ref{divgen}).
\vspace{0.5cm}
\section{Symmetry classes of tensors} \label{sec2}%
Let $\bbK$ be the field of real or complex numbers $\bbR$, $\bbC$.
We denote by $\bbK [{\calS}_r]$ the {\itshape group ring} of a symmetric group ${\calS}_r$. Furthermore we consider the $\bbK$-vector space ${\calT}_r V$ of $r$-times covariant $\bbK$-valued tensors $T$ over a finite dimensional $\bbK$-vector space $V$. Every group ring element $a = \sum_{p \in {\calS}_r} a(p)\,p \in \bbK [{\calS}_r]$ acts as so-called {\itshape symmetry operator} on tensors $T \in {\calT}_r V$ according to the definition
\begin{eqnarray}
(a T)(v_1 , \ldots , v_r) & := & \sum_{p \in {\calS}_r} a(p)\,
T(v_{p(1)}, \ldots , v_{p(r)}) \;\;\;\;\;,\;\;\;\;\;
v_i \in V \,. \label{symop}
\end{eqnarray}
Equation \eqref{symop} is equivalent to
\begin{eqnarray}
(a T)_{i_1 \ldots i_r} & = & \sum_{p \in {\calS}_r} a(p)\,
T_{i_{p(1)} \ldots  i_{p(r)}} \,.
\end{eqnarray}
\begin{Def}
Let $\frakr \subseteq \bbK [{\calS}_r]$ be a right ideal of $\bbK [{\calS}_r]$ for which an $a \in \frakr$ and a $T \in {\calT}_r V$ exist such that $aT \not= 0$. Then the tensor set
\begin{eqnarray}
{\calT}_{\frakr} & := & \{ a T \;|\; a \in \frakr \;,\;
T \in {\calT}_r V \}
\end{eqnarray}
is called the {\itshape symmetry class} of tensors defined by $\frakr$. If $\frakr$ is a minimal right ideal, then ${\calT}_{\frakr}$ is called {\it irreducible}.
\end{Def}
Since $\bbK [{\calS}_r]$ is semisimple for $\bbK = \bbR , \bbC$, every right ideal $\frakr \subseteq \bbK [{\calS}_r]$ possesses a generating idempotent $e$, i.e. $\frakr$ fulfils $\frakr = e \cdot \bbK [{\calS}_r]$.
\begin{Lem}\footnote{See H. Boerner \cite[p.127]{boerner} or B. Fiedler \cite{fie18}, \cite[p.110]{fie16}.}
If $e$ is a generating idempotent of $\frakr$, then a tensor $T \in {\calT}_r V$ belongs to ${\calT}_{\frakr}$ iff
\begin{eqnarray}
e T & = & T \,.
\end{eqnarray}
Thus we have
\begin{eqnarray}
{\calT}_{\frakr} & = & \{ eT \;|\; T \in {\calT}_r V \} \,. \label{classgen}
\end{eqnarray}
\end{Lem}
Symmetry classes can be characterized by left ideals of $\bbK [{\calS}_r]$, too. To see this, we construct group ring elements from tensors.
\begin{Def}
Every tensor $T \in {\calT}_r V$ and every $r$-tuple $b = (v_1 , \ldots , v_r) \in V^r$ of vectors from $V$ induce a group ring element
\begin{eqnarray}
T_b & := & \sum_{p \in {\calS}_r}\,T(v_{p(1)}, \ldots , v_{p(r)})\,p \;\;\in \bbK [{\calS}r]\,. \label{tb}
\end{eqnarray}
\end{Def}

A connection between (\ref{symop}) and (\ref{tb}) is given by the formula\footnote{See B. Fiedler \cite{fie18} or \cite[p.110]{fie16}.}
\begin{eqnarray}
\forall\,T \in {\calT}_r V \;,\; \forall\,a \in \bbK [{\calS}_r] \;,\; \forall\, b \in V^r \;: \;\;\;\; (aT)_b & = & T_b \cdot a^{\ast}
\end{eqnarray}
where the star '$\ast$' denotes the mapping
\begin{eqnarray}
\ast : a \;=\; \sum_{p \in {\calS}_r} a(p)\,p & \mapsto &
a^{\ast} \;:=\; \sum_{p \in {\calS}_r} a(p)\,p^{-1} \,.
\end{eqnarray}
Now, if a tensor $T$ belongs to a certain symmetry class, then its $T_b$ lie in a certain left ideal.
\begin{Prop}\footnote{See B. Fiedler \cite[Sec.III.3.1]{fie16} and \cite{fie8,fie14,fie17}.}
Let $e \in \bbK [{\calS}_r]$ be an idempotent. Then a
$T \in {\calT}_r V$ fulfils the condition $eT = T$ iff
$T_b \in \frakl := \bbK [{\calS}_r] \cdot e^{\ast}$ for all
$b \in V^r$, i.e. all $T_b$ of $T$ lie in the left ideal $\frakl$ generated by $e^{\ast}$. Moreover, if $\dim V \ge r$ then $\frakl$ is spanned by all group ring elements $T_b$ of tensors $T \in {\calT}_{{\frakl}^{\ast}}$, i.e.
$\frakl = \calL \{ T_b \;|\; T \in {\calT}_{{\frakl}^{\ast}} \,,\,
b \in V^r \}$ where $\calL$ denotes the linear closure.
\end{Prop}
Thus we can use the left ideal $\frakl = {\frakr}^{\ast}$ instead of the right ideal $\frakr$ to charakterize the symmetry class of $\frakr$.
\vspace{0.5cm}
\section{Symmetry classes of the first covariant derivatives of tensor fields} \label{sec3}%
Now we determine results about the symmetry classes of the first covariant derivatives of tensor fields. In particular, we are interested in symmetric or alternating covariant tensor fields. 

We consider only differentable objects of class $C^{\infty}$.
Let $M$ be an $m$-dimensional {\it differentable manifold} equipped with a {\it linear connection} or {\it covariant derivative} $\nabla$. We denote by $\calT_r M$, $r \ge 0$, the set of {\it differentable covariant tensor fields} of order $r$ on $M$.
If $T \in \calT_r M$, $r \ge 1$, then its {\it covariant derivative} $\nabla T$ has a coordinate representation\footnote{We use the Einstein summation convention, i.e. a symbol such as $T_{i_1 \ldots s \ldots\; \ldots i_r}^{\;\;\;\;\;\;\;\;\;\;s}$ means $T_{i_1 \ldots s \ldots\; \ldots i_r}^{\;\;\;\;\;\;\;\;\;\;s} := \sum_{s = 1}^m T_{i_1 \ldots s \ldots\; \ldots i_r}^{\;\;\;\;\;\;\;\;\;\;s}$.}
\begin{eqnarray}
{\nabla}_{i_{r+1}}T_{i_1 \ldots i_r} & = &
{\partial}_{i_{r+1}}T_{i_1 \ldots i_r} -
\sum_{k = 1}^r {\Gamma}_{i_{r+1} i_k}^{s_k}\,T_{i_1 \ldots i_{k-1} s_k i_{k+1} \ldots i_r} \label{covderiv1}
\end{eqnarray}
where $\partial$ is a partial derivative of the coordinates of $T$ and ${\Gamma}_{i j}^k$ are the connection coefficients of $\nabla$. Instead of (\ref{covderiv1}) we write also
\begin{eqnarray}
T_{i_1 \ldots i_r \,;\, i_{r+1}} & = &
T_{i_1 \ldots i_r \,,\, i_{r+1}} -
\sum_{k = 1}^r {\Gamma}_{i_{r+1} i_k}^{s_k}\,T_{i_1 \ldots i_{k-1} s_k i_{k+1} \ldots i_r}\,. \label{covderiv2}
\end{eqnarray}
Every tensor in a fixed point $p$ of $M$ can be gained as covariant derivative of a suitable tensor field.
\begin{Lem}
Let $M_p$ be the tangent space of $M$ in a point $p \in M$ of $M$ and
$W \in \calT_{r+1} M_p$, $r \ge 0$, be a covariant tensor of order $r + 1$ over $M_p$. Then we can find a covariant tensor field $T \in \calT_r M$ such that $(\nabla T)|_p = W$. \label{lem2}
\end{Lem}
\begin{proof}
First we consider the case $r \ge 1$. In a suitable open neighbourhood $U$ of $p$ we can choose a chart $x$ such that $x(p) = 0$. If $W_{i_1 \ldots i_{r+1}}$ are the coordinates of $W$ with respect to $x$ and $x^i$ are the coordinate functions of $x$ then $\tilde{T}_{i_1 \ldots i_r} := W_{i_1 \ldots i_r k} x^k$ yields a differentable tensor field on $U$. Further we can consider a function $\phi$ of class $C^{\infty}$ on $U$ for which open neighbourhoods $U_1$ and $U_2$ of $p$ exist such that $p \in U_1 \subset U_2 \subset U$ and $\phi |_{U_1} \equiv 1$, $\phi |_{U \setminus U_2} \equiv 0$. By means of $\phi$ we obtain a differentable tensor field $T \in \calT_r M$ if we set $T|_U := \phi \tilde{T}$ and $T|_{M \setminus U} := 0$. But $T$ fulfils $(\nabla T)|_p = W$ since we can write
\begin{eqnarray*}
T_{i_1 \ldots i_r \,;\, i_{r+1}}(p) & = &
T_{i_1 \ldots i_r \,,\, i_{r+1}}(p) -
\sum_{k = 1}^r {\Gamma}_{i_{r+1} i_k}^{s_k}(p)\,T_{i_1 \ldots i_{k-1} s_k i_{k+1} \ldots i_r}(p) \;=\; W_{i_i \ldots i_r i_{r+1}}\,.
\end{eqnarray*}
The last equality follows from $T_{i_1 \ldots  i_r}(p) = 0$ and $T_{i_1 \ldots i_r \,,\, i_{r+1}}(p) = W_{i_i \ldots i_r i_{r+1}}$.

In the case $r = 0$ the tensor $W$ has order 1 and $\tilde{T} = W_k x^k$ is a tensor field of order 0, i.e. a differentiable function. Obviously we can form the tensor field (function) $T \in \calT_0 M$ in the same way as in the case $r \ge 1$ and we obtain $T_{; i} = T_{, i} = W_i$ on the neighbourhood $U_1$ of $p$.
\end{proof}
Lemma \ref{lem2} leads to the consequence that every symmetry class can be generated by covariant derivatives of suitable tensor fields.
\begin{Cor}
Let $\frakr \subseteq \bbK [\calS_{r+1}]$, $r \ge 0$, be a right ideal with generating idempotent $e \in \bbK [\calS_{r+1}]$ for which an $a \in \frakr$ and a $W \in \calT_{r+1} M_p$ exist such that $aW \not= 0$. Then the symmetry class $\calT_{\frakr}$ of tensors of order $r +1$ over $M_p$ fulfils
\begin{eqnarray}
\calT_{\frakr} & = & \{ e (\nabla T)|_p \; | \; T \in \calT_r M \} \,.
\end{eqnarray} \label{cor1}
\end{Cor}
\begin{proof}
According to (\ref{classgen}) we have
$\calT_{\frakr} = \{ eW \;|\; W \in \calT_{r+1} M_p \}$. But for every $W \in \calT_{r+1} M_p$ there exists a $T \in \calT_r M$ such that $W = (\nabla T)|_p$. Thus Corollary \ref{cor1} follows.
\end{proof}
The {\it symmetrisation} '${\,}_{( \ldots )}$' and {\it anti-symmetrization} '${\,}_{[ \ldots ]}$' of a tensor field $T$ of order $r$ is defined by
\begin{eqnarray}
T_{(i_1 \ldots i_r)} & := & \frac{1}{r!}\,\sum_{p \in {\calS}_r}\,
T_{i_{p(1)} \ldots i_{p(r)}} \\
T_{[i_1 \ldots i_r]} & := & \frac{1}{r!}\,\sum_{p \in {\calS}_r}\,
\sign(p)\,T_{i_{p(1)} \ldots i_{p(r)}}\,.
\end{eqnarray}
From now on we consider only covariant derivatives $\nabla$ which are {\it torsion-free}, i.e. $\Gamma_{[i j]}^k = 0$. It is well-known for such $\nabla$ that the operations '${\,}_{( \ldots )}$' or '${\,}_{[ \ldots ]}$' and the operator $\nabla$ are permutable. However, this statement is correct for arbitrary symmetry operators $a \in \bbK [\calS_r]$, too.
\begin{Lem} Let $\nabla$ be torsion-free and $a = \sum_{p \in \calS_r} a(p)\,p \in \bbK [\calS_r]$, $r \ge 2$.
If we consider a $T \in \calT_r M$ and set $H_{i_1 \ldots i_r} := (aT)_{i_1 \ldots i_r}$  and
$W_{i_1 \ldots i_r i_{r+1}} := T_{i_1 \ldots i_r \,;\, i_{r+1}}$ with respect to arbitrary local coordinates, then it holds
\begin{eqnarray}
\sum_{p \in {\calS}_r}\,
a(p)\,W_{i_{p(1)} \ldots i_{p(r)} \, i_{r+1}} & = & H_{i_1 \ldots i_r \,;\, i_{r+1}}\,.
\end{eqnarray} \label{lem4}
\end{Lem}
\begin{proof}
Let $q \in M$ be an arbitrary point of $M$. We can choose such coordinates arround $q$ that all ${\Gamma}_{i j}^k$ vanish in $q$, i.e. ${\Gamma}_{i j}^k (q) = 0$. Then we have $W_{i_1 \ldots i_r i_{r+1}}(q) = T_{i_1 \ldots i_r \,,\, i_{r+1}}(q)$ and $H_{i_1 \ldots i_r \,;\, i_{r+1}}(q) = H_{i_1 \ldots i_r \,,\, i_{r+1}}(q)$. But since obviously
\begin{eqnarray*}
\sum_{p \in {\calS}_r}\,
a(p)\,T_{i_{p(1)} \ldots i_{p(r)} \,,\, i_{r+1}}(q) & = & H_{i_1 \ldots i_r \,,\, i_{r+1}}(q)\,,
\end{eqnarray*}
we obtain Lemma \ref{lem4}.
\end{proof}
The version of Lemma \ref{lem4} for the symmetry operators '${\,}_{( \ldots )}$' and '${\,}_{[ \ldots ]}$' reads
\begin{Lem} Let $\nabla$ be torsion-free and $r \ge 2$.
If we set $S_{i_1 \ldots i_r} := T_{(i_1 \ldots i_r)}$,
$A_{i_1 \ldots i_r} := T_{[i_1 \ldots i_r]}$ and
$W_{i_1 \ldots i_r i_{r+1}} := T_{i_1 \ldots i_r \,;\, i_{r+1}}$, then it holds
\begin{eqnarray}
W_{(i_1 \ldots i_r )\, i_{r+1}} & = & S_{i_1 \ldots i_r \,;\, i_{r+1}} \\
W_{[i_1 \ldots i_r ]\, i_{r+1}} & = & A_{i_1 \ldots i_r \,;\, i_{r+1}}\,.
\end{eqnarray} \label{lem1}
\end{Lem}
\begin{Prop} Consider the case $r \ge 2$.
Let $\tilde{\calS}_r := \{ p \in \calS_{r+1} \;|\; p(r+1) = r+1 \}$ be the subgroup of those permutations of $\calS_{r+1}$ which have $r+1$ as fixed point. Then
\begin{eqnarray}
e_s \; := \; \frac{1}{r!}\,\sum_{p \in \tilde{\calS}_r} \,p
& \;\;\;,\;\;\; &
e_a \; := \; \frac{1}{r!}\,\sum_{p \in \tilde{\calS}_r} \,\sign(p)\,p
\end{eqnarray}
are idempotents of the group ring $\bbK [\calS_{r+1}]$ which fulfil
\begin{eqnarray}
e_s^{\ast} \; = \; e_s & \;\;\;,\;\;\; & e_a^{\ast} \; = \; e_a \,.
\end{eqnarray}
If $\nabla$ is a torsion-free covariant derivative and $S$, $A$ are symmetric or alternating differentiable tensor fields of order $r$, respectively, then it holds\footnote{Here we assume that the tensor indices are numbered in the manner of (\ref{covderiv2}).}
\begin{eqnarray}
e_s^{\ast} \nabla S \; = \; \nabla S & \;\;\;,\;\;\; & e_a^{\ast} \nabla A \; = \; \nabla A \,. \label{eee}
\end{eqnarray} \label{prop1}
\end{Prop}
\begin{proof}
Taking into account Lemma \ref{lem1} and
$\mathrm{sign}(p \cdot q) = \mathrm{sign}(p) \mathrm{sign}(q)$,
$\mathrm{sign}(p^{-1}) = \mathrm{sign}(p)$, we can prove Proposition \ref{prop1} by simple calculations.
\end{proof}
Now we will show that the symmetry classes of covariant derivatives of tensor fields are characterized by {\it Littlewood-Richardson products}.
Let us denote by $\regrep$ the {\it regular representation} of the symmetric group $\calS_{r+1}$:
\begin{eqnarray}
\regrep: \calS_{r+1} \rightarrow \mathrm{Gl}(\bbK [\calS_{r+1}])
& \;,\; &
\regrep_p (a) \,=\, p \cdot a \;,\;
p \in \calS_{r+1} \;,\; a \in \bbK [\calS_{r+1}]\,.
\end{eqnarray}
The left ideals $\frakl \subseteq \bbK [\calS_{r+1}]$ of $\bbK [\calS_{r+1}]$ can be considered representation spaces of subrepresentations $\rho = \regrep |_{\frakl}$ of $\regrep$:
\begin{eqnarray}
\rho: \calS_{r+1} \rightarrow \mathrm{Gl}(\frakl)
& \;,\; &
{\rho}_p (a) \,=\, p \cdot a \;,\;
p \in \calS_{r+1} \;,\; a \in \frakl\,.
\end{eqnarray}
We see from a generalization of Proposition \ref{prop1} that investigations of covariant derivatives of tensor fields can be based on following
\begin{Set}
Let $T \in \calT_r M$ be a differentiable tensor field of order $r$ on $M$ whose symmetry class is defined by a left ideal $\bbK [\calS_{r}] \cdot e$ with generating idempotent $e \in \bbK [\calS_{r}]$, i.e. $e^{\ast} T = T$. We identify $\calS_r$ with $\tilde{\calS}_r$ by means of $[i_1 , \ldots , i_r] \mapsto [i_1 , \ldots , i_r , r+1]$ (where $[i_1 , \ldots , i_r]$ is the list representation of a permutation) and denote by $\tilde{e} \in \bbK[\tilde{\calS}_r]$ the corresponding embedding of $e \in \bbK [\calS_{r}]$ into $\bbK [\calS_{r+1}]$. If $\nabla$ is a torsion-free covariant derivative on $M$ then the left ideal $\frakl := \bbK [\calS_{r+1}] \cdot \tilde{e}$ defines the symmetry class of $\nabla T$ because $\tilde{e}^{\ast} \nabla T = \nabla T$. Furthermore, we can consider the left ideal $\tilde{\frakl} := \bbK [\tilde{\calS}_{r}] \cdot \tilde{e}$ of $\bbK [\tilde{\calS}_{r}]$. The left ideals $\frakl$ and $\tilde{\frakl}$ are the representation spaces of the representations\footnote{If $\alpha$ is a representation of a group $G$ and $H \subseteq G$ is a subgroup of $G$, then $\alpha \downarrow H$ denotes the restriction of $\alpha$ to $H$.}
\begin{eqnarray}
\frakl & \Longleftrightarrow & \rho \;:=\; \regrep |_{\frakl} \\
\tilde{\frakl} & \Longleftrightarrow & \sigma \;:=\; (\regrep  \downarrow \tilde{\calS}_r) |_{\tilde{\frakl}} \,.
\end{eqnarray} \label{setting}
\end{Set}
A simple consequence of Corollary \ref{cor1}, Lemma \ref{lem4}  and Setting \ref{setting} is that the symmetry class of the above $\nabla T$ is generated by the covariant derivatives of the symmetizations $e^{\ast} W$ of arbitrary tensor fields $W \in \calT_r M$.
\begin{Prop}
Assume that Setting {\rm \ref{setting}} and $r \ge 1$ are valid and $p \in M$. Then we have
\begin{eqnarray*}
\calT_{\frakl^{\ast}} & = &
\left\{ \tilde{e}^{\ast}(\nabla W)|_p \;|\; W \in \calT_r M \right\}
\;=\;
\left\{ (\nabla e^{\ast} W)|_p \;|\; W \in \calT_r M \right\}\,.
\end{eqnarray*} \label{prop2}
\end{Prop}
Now we determine Littlewood-Richardson products describing covariant derivatives.
\begin{Thm}
Assume that Setting {\rm \ref{setting}} and $r \ge 1$ are valid. Then the representation $\rho$ is a Littlewood-Richardson product
\begin{eqnarray}
\rho & \sim & \sigma \, [1] \,.
\end{eqnarray}
If $\sigma$ is an irreducible representation $\sigma \sim [\lambda]$, $\lambda \vdash r$, then the Littlewood-Richardson rule\footnote{See
D. E. Littlewood \cite[pp.94-96]{littlew1},
A. Kerber \cite[Vol.240/p.84]{kerber},
G. D. James and A. Kerber \cite[p.93]{jameskerb},
A. Kerber \cite[Sec.5.5]{kerber3},
I. G. Macdonald \cite[Chap.I,Sec.9]{mcdonald},
R. Merris \cite[p.100]{merris},
W. Fulton and J. Harris \cite[pp.455-456]{fultharr},
S. A. Fulling, R. C. King, B. G. Wybourne and C. J. Cummins \cite{full4}. See also B. Fiedler \cite[Sec.II.5]{fie16}.}
 yields a multiplicity-free decomposition
\begin{eqnarray}
\rho \;\sim\; [\lambda][1] & \sim & \sum_{\twoline{\mu \vdash r+1}{\lambda \subset \mu}} \,[\mu]\,, \label{lrsum}
\end{eqnarray}
i.e. $[\lambda][1]$ decomposes into a sum of all Young frames $[\mu]$ which can be formed from $[\lambda]$ by adding one box to $[\lambda]$.
\label{thm2}
\end{Thm}
\begin{proof}
When we introduce the notation $\tilde{\calS}_1 := \{ \id \}$ for the trivial subgroup of $\calS_{r+1}$, then the set product $\tilde{\calS_r} = \tilde{\calS_r} \cdot \tilde{\calS}_1$ is a direct product $\tilde{\calS_r} = \tilde{\calS_r} \times \tilde{\calS}_1$. (It is even a Young subgroup.) We consider the representations
\begin{eqnarray}
\iota: \tilde{\calS}_1 \rightarrow \mathrm{Gl}(\bbK [\tilde{\calS}_1])
& \;,\; &
{\iota}_{\id}(u) = u \;\;,\;\; u \in \bbK [\tilde{\calS}_1] \\
\sigma: \tilde{\calS}_r \rightarrow \mathrm{Gl}(\bbK [\tilde{\calS}_r] \cdot \tilde{e})
& \;,\; &
{\sigma}_p (v) = p \cdot v \;\;,\;\; p \in \tilde{\calS}_r \;,\; v \in \bbK [\tilde{\calS}_r] \cdot \tilde{e} \,.
\end{eqnarray}
Obviously we can regard $\sigma$ as an outer tensor product of representations
\begin{eqnarray*}
\sigma & = & \sigma \,\#\, \iota \,.
\end{eqnarray*}
Thus the left ideal $\frakl = \bbK [\calS_{r+1}] \cdot \tilde{e}$ is the representation space of the induced representation
$(\sigma \,\#\, \iota) \uparrow \calS_{r+1}$  which has the structure of a Littlewood-Richardson product, i.e.
\begin{eqnarray*}
\rho & = & (\sigma \,\#\, \iota) \uparrow \calS_{r+1} \;\sim\; \sigma [1] \,.
\end{eqnarray*}
If $\sigma$ is irreducible, i.e. $\sigma \sim [\lambda]$, $\lambda \vdash r$, then the Littlewood-Richardson rule yields (\ref{lrsum})
.
\end{proof}
\begin{Rem}
If we restrict us to  irreducible representations $\sigma$, then the proof of Theorem \ref{thm2} is a repetition of a part of the proof of the {\it branching theorem} for irreducible representations of symmetric groups (see A. Kerber \cite[Vol.240/p.85]{kerber}).
\end{Rem}
\begin{Rem}
If $\sigma$ is a reducible representation and we know a decomposition
$\sigma = \bigoplus_i {\sigma}_i$ into subrepresentations ${\sigma}_i$ (irreducible or reducible), then we can use the formula
\begin{eqnarray}
\sigma [1] & \sim & \sum_i \, {\sigma}_i [1]
\end{eqnarray}
to determine the structure of the decomposition of $\sigma [1]$ into irreducible subrepresentations.
\end{Rem}
\begin{Rem}
Additional information about a tensor field considered can lead to a further reduction of the sum (\ref{lrsum}). For instance it is well-known that the symmetry classes of the {\it Riemannian curvature tensor} $R$ and its covariant derivative\footnote{Here we assume that $\nabla$ is the Levi-Civita connection of a pseudi-Riemannian fundamental tensor $g \in \calT_2 M$ and $R$ is the curvature tensor of $\nabla$.} $\nabla R$ are defined by the {\it Young symmetrizers}\footnote{See Section \ref{sec4} for some details about Young symmetrizers.} $y_t$ and $y_{t'}$ of the Young tableaux\footnote{See S. A. Fulling, R. C. King, B. G. Wybourne and C. J. Cummins \cite{full4}. See also B. Fiedler \cite{fie12}.}
\begin{eqnarray}
t \;=\;
\begin{array}{|c|c|}
\hline
1 & 3 \\
\hline
2 & 4 \\
\hline
\end{array}
&\;\;\;,\;\;\; &
t' \;=\;
\begin{array}{|c|c|c|c}
\cline{1-3}
1 & 3 & 5 & \\
\cline{1-3}
2 & 4 &\multicolumn{2}{c}{\;\;\;} \\
\cline{1-2}
\end{array}
\,. \label{riemtabls}
\end{eqnarray}
However, if we apply (\ref{lrsum}) to the tableau $t$ we obtain
$[2^2][1] \sim [3,2] + [2^2,1] \not\sim [3,2]$.
The difference results from the fact that $\nabla R$ fulfils the second Bianchi identity
\begin{eqnarray*}
R_{i j k l \,;\, m} \;+\; R_{i j l m \,;\, k} \;+\; R_{i j m k \,;\, l} & = & 0
\end{eqnarray*}
which is not satisfied by other tensor fields from $\calT_4 M$ in general.

A second example which shows such effects is the case of higher covariant derivatives of tensor fields. If we apply Theorem \ref{thm2} to covariant derivatives\footnote{Again we assume that $\nabla$ is a Levi-Civita connection.} of second order $T_{i_1 \ldots i_r \,;\, i_{r+1} i_{r+2}}$ of a tensor field $T \in \calT_r M$ then Theorem \ref{thm2} yields a result in which the so-called {\it Ricci identity}
\begin{eqnarray}
T_{i_1 \ldots i_r \,;\, [i_{r+1} i_{r+2}]} & = &
\frac{1}{2}\, \sum_{k = 1}^r \, R_{i_{r+1} i_{r+2} i_k}^{\hspace*{40pt} s_k}\,
T_{i_1 \ldots i_{k-1} s_k i_{k+1} \ldots i_r} \label{ricciid}
\end{eqnarray}
was left out of account. Thus the set of Young frames determined by multiple application of (\ref{lrsum}) will be ''too large''. (\ref{lrsum}) produces a set of Young frames which is correct also for covariant derivatives
$T_{i_1 \ldots i_r i_{r+1} \,;\, i_{r+2}}$ of tensor fields $T \in \calT_{r+1} M$ of order $r + 1$ to which an identity (\ref{ricciid}) is irrelevant.
\end{Rem}
Now we present a version of Theorem \ref{thm2} for the special case of symmetric or alternating tensor fields.
\begin{Thm}
Assume that $r\ge 2$. Let $\frakl_s := \bbK [\calS_{r+1}] \cdot e_s$,
$\frakl_a := \bbK [\calS_{r+1}] \cdot e_a$ be the left ideals generated by the idempotents $e_s$, $e_a$. Then the subrepresentations
${\rho}_s := \regrep |_{\frakl_s}$, ${\rho}_a := \regrep |_{\frakl_a}$ are Littlewood-Richardson products ${\rho}_s \sim [r][1]$, ${\rho}_a \sim [1^r][1]$ for which the Littlewood-Richardson rule yields decompositions
\begin{eqnarray}
[r][1] \;\sim\; [r+1] + [r , 1]
& \;\;\;,\;\;\; &
[1^r][1] \;\sim\; [1^{r+1}] + [2 , 1^{r-1}]\,. \label{lrrules}
\end{eqnarray}
The idempotents $e_s$, $e_a$ have unique decompositions corresponding to {\rm (\ref{lrrules})} into primitive orthogonal idempotents
\begin{eqnarray}
e_s \;=\; f_s + h_s & \;\;\;,\;\;\; e_a \;=\; f_a + h_a \label{idpotdec}
\end{eqnarray}
which generate the minimal left ideals that are the representation spaces of the irreducible representations in {\rm (\ref{lrrules})}:
\begin{eqnarray*}
\bbK [\calS_{r+1}] \cdot f_s & \Longleftrightarrow & [r+1]\\
\bbK [\calS_{r+1}] \cdot h_s & \Longleftrightarrow & [r , 1]\\
\bbK [\calS_{r+1}] \cdot f_a & \Longleftrightarrow & [1^{r+1}]\\
\bbK [\calS_{r+1}] \cdot h_a & \Longleftrightarrow & [2 , 1^{r-1}]\,.
\end{eqnarray*}
We know all idempotents in {\rm (\ref{idpotdec})} since
\begin{eqnarray}
f_s \; = \; \frac{1}{(r+1)!}\,\sum_{p \in \calS_{r+1}} \,p
& \;\;\;,\;\;\; &
f_a \; = \; \frac{1}{(r+1)!}\,\sum_{p \in \calS_{r+1}} \,\sign(p)\,p \,.
\label{idpotf}
\end{eqnarray} \label{thm1}
\end{Thm}
\begin{proof}
If we apply Theorem \ref{thm2} to the representations
\begin{eqnarray}
{\sigma}_s: \tilde{\calS}_r \rightarrow \mathrm{Gl}(\bbK [\tilde{\calS}_r] \cdot e_s)
& , &
({\sigma}_s)_p (v) = p \cdot v \;\;,\;\; p \in \tilde{\calS}_r \;,\; v \in \bbK [\tilde{\calS}_r] \cdot e_s \\
{\sigma}_a: \tilde{\calS}_r \rightarrow \mathrm{Gl}(\bbK [\tilde{\calS}_r] \cdot e_a)
& , &
({\sigma}_a)_p (v) = p \cdot v \;\;,\;\; p \in \tilde{\calS}_r \;,\; v \in \bbK [\tilde{\calS}_r] \cdot e_a
\end{eqnarray}
then we obtain (\ref{lrrules}) by means of the {\it Littlewood-Richardson rule} since ${\sigma}_s$ and ${\sigma}_a$ are irreducible and fulfil ${\sigma}_s \sim [r]$, ${\sigma}_a \sim [1^r]$.

The relations (\ref{lrrules}) tell us that every of the representations $({\sigma}_s \,\#\, \iota) \uparrow \calS_{r+1}$, $({\sigma}_a \,\#\, \iota) \uparrow \calS_{r+1}$ decomposes into two irreducible representations which have multiplicities 1. Thus the generating idempotents $e_s$, $e_a$ of the representation spaces possess unique corresponding decompositions (\ref{idpotdec}) into primitive othogonal idempotents. Since we know that the idempotents for representations $[r+1]$ and $[1^{r+1}]$ are given by (\ref{idpotf}) we can calculate the remaining idempotents $h_s$ and $h_a$, too.
\end{proof}
\begin{Rem}
In the case of an alternating tensor field $A \in \calT_r M$ the symmetry operator $f_a$ transforms $\nabla A$ into the exterior derivative $dA$ of $A$, i.e.
$f_a^{\ast} (\nabla A) = f_a (\nabla A) = dA$. Thus the symmetry operator $h_a$ yields the difference of $\nabla A$ and $dA$, i.e. $h_a^{\ast} (\nabla A) = \nabla A - dA$.
\end{Rem}
Now let us consider the more general case of an arbitrary tensor field $T \in \calT_r M$ whose symmetry class is defined by a known primitive idempotent $e \in \bbK [\calS_r]$.
Also in this case, there is a simple possibility to calculate all primitive idempotents which belong to a decomposition (\ref{lrsum}) for the covariant derivatives of $T$. A starting point is the well-known
\begin{Lem}\footnote{See H. Boerner \cite[Sec.III.3, III.4]{boerner} and R. Merris \cite[Sec.4]{merris}, in particular Exercise 41 in \cite[p.117]{merris}. See also B. Fiedler \cite[Prop.II.1.47, Prop.I.1.8]{fie16}.}
Let $\lambda \vdash r$ be a partition of $r \ge 1$ and ${\chi}_{\lambda}$ be the irreducible character of $\calS_r$ which belongs to that equivalence class of irreducible representations of $\calS_r$ which contains the irreducible representations\footnote{Here $\regrep$ denotes the regular representation of $\calS_r$.} $\regrep |_{\bbK [\calS_r] \cdot y_t}$ defined by Young symmetrizers $y_t$ with Young frame $\lambda$. Then the group ring element
\begin{eqnarray}
e_{\lambda} & := & \frac{{\chi}_{\lambda}(id)}{r!}\,\sum_{p \in \calS_r} {\chi}_{\lambda} (p)\,p^{-1}
\end{eqnarray}
is the unique centrally primitive idempotent that generates the minimal two-sided ideal $\fraka_{\lambda} := \bbK [\calS_r]\cdot e_{\lambda}$ from the isotypic decomposition $\bbK [\calS_r] = \bigoplus_{\mu \vdash r} \fraka_{\mu}$ of the group ring $\bbK [\calS_r]$ into minimal two-sided ideals $\fraka_{\mu}$. $\fraka_{\lambda}$ is that minimal two-sided ideal which contains all left ideals $\bbK [\calS_r] \cdot y_t$ generated by Young symmetrizers $y_t$ with Young frame $\lambda$. The idempotents $e_{\lambda}$ fulfil
\begin{eqnarray}
\sum_{\lambda \vdash r}\,e_{\lambda} \;=\; \id
& \;\;\;\;\; ; \;\;\;\;\;
e_{\lambda} \cdot e_{\lambda'} \;=\; 0 \;\;\;\mathrm{if}\;\; \lambda \not= \lambda'\,.
\end{eqnarray}
\end{Lem}
\begin{Thm}
Assume, that the symmetry class of a tensor field $T \in \calT_r M$, $r \ge 1$, is defined by a primitive idempotent $e \in \bbK [\calS_r]$ whose representation $\sigma$ according to Setting \ref{setting} satisfies $\sigma \sim [\lambda]$, $\lambda \vdash r$. Then \begin{eqnarray}
\tilde{e} \;=\; \sum_{\twoline{\mu \vdash r + 1}{\lambda \subset \mu}}\,h_{\mu}
& \;\;\;,\;\;\; &
h_{\mu} \;:=\; \tilde{e} \cdot e_{\mu} \label{decomp1}
\end{eqnarray}
yields the decomposition of $\tilde{e}$ into primitive idempotents $h_{\mu}$, which corresponds to relation {\rm (\ref{lrsum})}.
\end{Thm}
\begin{proof}
Because (\ref{lrsum}) is multiplicity-free, a decomposition of $\tilde{e}$ according to (\ref{lrsum}) into primitive idempotents contains exactly one primitive idempotent $h_{\mu}$ for every $[\mu]$ in (\ref{lrsum}). Every such $h_{\mu}$ lies in the corresponding two-sided ideal $\fraka_{\mu}$, i.e. $h_{\mu} \in \fraka_{\mu}$. On the other hand, we can write
$\tilde{e} = \tilde{e} \cdot \id = \sum_{\mu \vdash r+1} \tilde{e} \cdot e_{\mu}$. Since $\tilde{e} \cdot e_{\mu} \in \fraka_{\mu}$, we obtain (\ref{decomp1}).
\end{proof}
If we carry out calculations in large $\calS_r$, then a use of formula (\ref{decomp1}) leads to very high costs in calculation time and computer memory. However, fast discrete Fourier transforms can help to solve this problem.
\begin{Def}
A {\it discrete Fourier transform} for $\calS_r$ is an isomorphism
\begin{eqnarray}
D : \; \bbK [\calS_r] & \rightarrow &
\bigotimes_{\lambda \vdash r} {\bbK}^{d_{\lambda} \times d_{\lambda}}
\end{eqnarray}
according to Wedderburn's theorem which maps the group ring $\bbK [\calS_r]$ onto an outer direct product $\bigotimes_{\lambda \vdash r} {\bbK}^{d_{\lambda} \times d_{\lambda}}$ of full matrix rings ${\bbK}^{d_{\lambda} \times d_{\lambda}}$. We denote by $D_{\lambda}$ the {\it natural projections}
$D_{\lambda} : \bbK [\calS_r] \rightarrow
{\bbK}^{d_{\lambda} \times d_{\lambda}}$.
\end{Def}
Since the subrings $(0 , \ldots , 0 , {\bbK}^{d_{\lambda} \times d_{\lambda}} ,
0 , \ldots , 0 )$ correspond to the two-sided ideals $\fraka_{\lambda}$ under a discrete Fourier transform we obtain
\begin{Cor}
If a discrete Fourier transform $D$ is known for $\bbK [\calS_{r+1}]$, then the idempotents $h_{\mu}$ in {\rm (\ref{decomp1})} can be calculated by
\begin{eqnarray}
h_{\mu} & = & D^{-1} \left( (0,\ldots , D_{\mu}(\tilde{e}) , \ldots , 0 ) \right) \,. \label{fft}
\end{eqnarray}
\end{Cor}
A use of (\ref{fft}) in computer calculations is much more efficient than an application of (\ref{decomp1}). A very efficient algorithm of a fast Fourier transform for $\calS_r$ was developped by M. Clausen and U. Baum (see \cite{clausbaum1,clausbaum2}). It is based on {\it Young's seminormal representation} of $\calS_r$. Our {\sf Mathematica} package {\sf PERMS} \cite{fie10} uses {\it Young's natural representation} of $\calS_r$ as discrete Fourier transform.
\vspace{0.5cm}
\section{Do Young symmetrizers describe the symmetry classes of $\nabla S$ or $\nabla A$?} \label{sec4}
Now we turn to the question wheter the symmetry classes of $\nabla S$ and $\nabla A$ can be characterized by {\it Young symmetrizers}. According to Proposition \ref{prop1} and Theorem \ref{thm1} the symmetry classes of $\nabla S$ and $\nabla A$ are defined by the idempotents
\begin{eqnarray*}
e_s \;=\; f_s + h_s & \;\;\;,\;\;\; e_a \;=\; f_a + h_a\,.
\end{eqnarray*}
The idempotents $f_s$ and $f_a$ are proportional to the Young symmetrizers of the Young frames $(r + 1) \vdash r + 1$ or $(1^{r + 1}) \vdash r + 1$, respectively. Now we investigate the
\begin{Pro}
Can we find Young tableaux $t_s$, $t_a$ with frame $(r , 1) \vdash r + 1$ or $(2 , 1^{r-1}) \vdash r + 1$, respectively, such that the idempotents $e_{t_s} = {\mu}_{t_s} y_{t_s}$, $e_{t_a} = {\mu}_{t_a} y_{t_a}$ are generating idempotents of the minimal left ideals $\frakl'_s := \bbK [\calS_{r+1}] \cdot h_s$,
$\frakl'_a := \bbK [\calS_{r+1}] \cdot h_a$? Here $y_{t_s}$ and $y_{t_a}$ are the {\it Young symmetrizers}\footnote{We define a Young symmetrizer $y_t$ of a Young tableau $t$ by the formula $y_t := \sum_{p \in \calH_t} \sum_{q \in \calV_t} \mathrm{sign}(q)\,p \cdot q$, where $\calH_t$, $\calV_t$ are the groups of horizontal or vertical permutations of $t$, respectively.} of the Young tableaux $t_s$, $t_a$ and ${\mu}_{t_s}, {\mu}_{t_a} \not= 0$ are constants.
\end{Pro}
If the answer is ''yes'', then the symmetry classes of $\nabla S$ or $\nabla A$ are determined by new idempotents
\begin{eqnarray}
\tilde{e}_s \;=\; f_s + e_{t_s} & \;\;\;,\;\;\; \tilde{e}_a \;=\; f_a + e_{t_a}\,,
\end{eqnarray}
which are completely built from Young symmetrizers.

We investigated this problem by computer calculations by means of the {\sf Mathematica} package {\sf PERMS} \cite{fie10} in the groups $\calS_3$, $\calS_4$, $\calS_5$, i.e.in the cases $r = 2, 3, 4$. We obtained the following results:
\begin{itemize}
\item{No Young symmetrizer idempotent $e_t$ of a Young frame $(r , 1) \vdash r + 1$ reproduces or annihilates $h_s \in \bbK [\calS_{r+1}]$, i.e. no relation
$h_s \cdot e_t = h_s$ or $h_s \cdot e_t = 0$ is satisfied.}
\item{For $h_a \in \bbK [\calS_{r+1}]$ there are many Young symmetrizer idempotents $e_t$ with Young frame $(2 , 1^{r-1}) \vdash r + 1$ which reproduce or annihilate $h_a$, i.e. which fulfil $h_a \cdot e_t = h_a$ or $h_a \cdot e_t = 0$. In particular, the idempotent $e_t$ of the lexicographically greatest standard tableau
\begin{eqnarray}
t & = &
\begin{array}{|c|c|c}
\cline{1-2}
1 & r + 1 & \\
\cline{1-2}
2 & \multicolumn{2}{c}{\;\;\;} \\
\cline{1-1}
\vdots & \multicolumn{2}{c}{\;\;\;} \\
\cline{1-1}
r & \multicolumn{2}{c}{\;\;\;} \\
\cline{1-1}
\end{array} \label{reptabl}
\end{eqnarray}
reproduces $h_a$ whereas the idempotents $e_t$ of all other standard tableaux of $(2 , 1^{r-1})$ annihilate $h_a$.}
\end{itemize}
For the $\calS_3$ and  $r = 2$, we also verified these results by calculations by means of the packages {\sf PERMS} \cite{fie10} and {\sf Ricci} \cite{ricci3} in which we checked the action of idempotents $e_t^{\ast}$ onto tensors with a symmetry given by $h_s$ or $h_a$. {\sf Mathematica} notebooks of all above calculations can be downloaded from my internet page \cite{fie21}.
We present tables of all Young tableaux, whose idempotents $e_t$ reproduce or annihilate $h_a$, in the Appendix.

Now we present theorems which tell us that essential parts of the above computer results are valid for all $r \ge 2$.

It is well-known that two idempotents $e , f \in \bbK [\calS_r]$ generate the same left ideal iff $e \cdot f = e$ and $ f \cdot e = f$. In the case of primitive idempotents $e , f$ we have
\begin{Lem}
If $e , f \in \bbK [\calS_r]$ are primitive idempotents, then the equations
$e \cdot f = e$ and $f \cdot e = f$ are equivalent. \label{lem3}
\end{Lem}
\begin{proof}
Assume that $e \cdot f = e$ is valid. Then $e$ is an element of the left ideal
$\frakl = \bbK [\calS_r] \cdot f$ and generates a non-vanishing subideal $\frakl' = \bbK [\calS_r] \cdot e$ of $\frakl$. Since $e$ and $f$ are primitive, the left ideals $\frakl$ and $\frakl'$ are minimal. Thus we obtain $\frakl = \frakl'$. But then, it follows $f \cdot e = f$ because $f$ belongs to the left ideal $\frakl'$ generated by $e$.
\end{proof}
\begin{Thm}
Consider the idempotent $h_a$ for an arbitrary $r \ge 2$. Then the Young symmetrizer idempotent $e_t$ of the lexicographically greatest standard tableau {\rm (\ref{reptabl})} with Young frame $(2 , 1^{r-1}) \vdash r + 1$ reproduces $h_a$, i.e. $h_a \cdot e_t = h_a$. Furthermore the idempotents $e_t$ of all other standard tableaux $t$ with Young frame $(2 , 1^{r-1}) \vdash r + 1$ annihilate $h_a$, i.e. $h_a \cdot e_t = 0$.
\end{Thm}
\begin{proof}
The groups $\calH_t$ and $\calV_t$ of horizontal/vertical permutations of the Young tableau $t$ according to (\ref{reptabl}) fulfil $\calH_t = \langle (1\,,\,r+1) \rangle$ and
$\calV_t = \tilde{\calS}_r$. Thus we can write
\begin{eqnarray}
e_t & = & {\mu}_t\,\sum_{p \in \calH_t} \sum_{q \in \calV_t} \mathrm{sign}(q)\,p \cdot q \;=\; {\mu}_t\cdot r!\;\{\id + (1 \,,\, r+1)\}\cdot e_a\,. \label{et}
\end{eqnarray}
It holds\footnote{See e.g. H. Boerner \cite[p.98]{boerner} or W. M\"uller \cite[p.73]{muell}. See also B. Fiedler \cite[Sec.I.3.1]{fie16}.} $e_t \cdot f_a = 0$ since $e_t$ and $f_a$ are proportional to Young symmetrizers of the different Young frames $(2 , 1^{r-1})$ and $(1^{r+1})$. From this and (\ref{et}) we obtain
\begin{eqnarray*}
e_t \cdot h_a & = & e_t \cdot (e_a - f_a) \;=\; e_t \cdot e_a \;=\; e_t\,.
\end{eqnarray*}
Because $e_t$ and $h_a$ are primitive idempotents, Lemma \ref{lem3} yields 
$h_a \cdot e_t = h_a$.

It is well-known\footnote{See A. Kerber \cite[Vol.240/p.73]{kerber} or H. Boerner \cite[p.101]{boerner}. See also B. Fiedler \cite[Sec.I.3.1]{fie16}.}: If $y_{t_1}$ and $y_{t_2}$ are Young symmetrizers of two standard tableaux $t_1$, $t_2$ which possess the same Young frame, and $t_1$ is lexicographically smaller\footnote{A tableau $t_2$ is regarded as greater than a tableau $t_1$ (of the same Young frame), if the simultaneous run through the rows of both tableaux from left to right and from top to bottom reaches earlier in $t_2$ a number which is greater than the number on the corresponding place in $t_1$.} than $t_2$, then they satisfy $y_{t_2} \cdot y_{t_1} = 0$. Because (\ref{et}) is built from the lexicographically greatest standard tableau $t$ of $(2 , 1^{r-1})$, we obtain
\begin{eqnarray*}
h_a \cdot e_{t'} & = & h_a \cdot e_t \cdot e_{t'} \;=\; 0
\end{eqnarray*}
for the idempotent $e_{t'}$ of every other standard tableau of $(2 , 1^{r-1})$
\end{proof}
\begin{Thm}
Consider the idempotent $h_s$ for an arbitrary $r \ge 2$. Then it holds $h_s \cdot e_t \not= h_s$ for all Young tableaux $t$ with Young frame $(r , 1) \vdash r + 1$.
\end{Thm}
\begin{proof}
A Young tableau $t$ with a Young frame $(r , 1) \vdash r + 1$ has a form
\begin{eqnarray}
t & = &
\begin{array}{|c|c|c|c|c|c}
\cline{1-5}
k & \star & \star & \ldots & \star & \\
\cline{1-5}
l & \multicolumn{5}{c}{\;\;\;} \\
\cline{1-1}
\end{array}\,. \label{tabl2}
\end{eqnarray}
If we assume that the first column of such a tableau contains the numbers $k$ and $l$ as in (\ref{tabl2}), then the groups of horizontal or vertical permutations of $t$ read
\begin{eqnarray}
\calH_t & = & \left\{ p \in \calS_{r+1} \;|\; p(l) = l \right\} \;=:\; (\calS_{r+1})_l \\
\calV_t & = & \langle (k\,,\,l) \rangle \,.
\end{eqnarray}

First we consider the case, that the first column of (\ref{tabl2}) does not contain $r+1$. In this case we have
\begin{eqnarray*}
y_t & = & \Bigl( \sum_{p \in (\calS_{r+1})_l} \,p \Bigr) \cdot \{ \id - (k\,,\,l) \}\,.
\end{eqnarray*}
But because $\calV_t \subseteq \tilde{\calS}_r$, we obtain $\{ \id - (k\,,\,l) \} \cdot e_s = 0$ and $y_t \cdot e_s = 0$. Thus $e_t = {\mu}_t y_t$ does not lie in the left ideal $\frakl_s := \bbK [\calS_{r+1}] \cdot e_s$ and can not play the role of a generating idempotent of $\frakl'_s := \bbK [\calS_{r+1}] \cdot h_s$. Consequently $h_s \cdot e_t \not= h_s$.

Netxt we investigate the case $l = r+1$. In this case we have
$\calH_t = \tilde{\calS}_r$, $\calV_t = \langle (k\,,\,r+1) \rangle$,
\begin{eqnarray}
e_t & = & {\mu}_t\,\Bigl( \sum_{p \in \tilde{\calS}_r} \,p \Bigr) \cdot \{ \id - (k\,,\,r+1) \} \;=\; {\mu}_t\cdot r!\,e_s \cdot \{ \id - (k\,,\,r+1) \}\,, \label{et2}
\end{eqnarray}
which leads to
\begin{eqnarray*}
h_s \cdot e_t & = & (e_s - f_s) \cdot e_t \;=\; e_s \cdot e_t \;=\; e_t \,.
\end{eqnarray*}
We decompose $h_s$ and $e_t$ into parts that correspond to the right cosets of $\calS_{r+1}$ relative to $\tilde{\calS}_r$. Obviously
$\fR := \{ (i\,,\,r+1) \;|\; i = 1 , \ldots , r+1 \}$ is a complete set of representatives of those right cosets. If we arrange the summands of $f_s$ and $h_s$ according to the decomposition of $\calS_{r+1}$ into cosets we obtain
\begin{eqnarray*}
f_s & = & \frac{1}{(r+1)!}\,\sum_{s \in \fR} \Bigl(
\sum_{p \in \tilde{\calS}_r}\,p \Bigr) \cdot s
\end{eqnarray*}
and
\begin{eqnarray*}
h_s \;=\; e_s - f_s & = & \frac{r}{(r+1)!}\,\Bigl(
\sum_{p \in \tilde{\calS}_r}\,p \Bigr)  \;-\; \frac{1}{(r+1)!}\,\sum_{\twoline{s \in \fR}{s \not= \id}} \Bigl(
\sum_{p \in \tilde{\calS}_r}\,p \Bigr) \cdot s\,.
\end{eqnarray*}
Thus $h_s$ has a non-vanishing part in every right coset of $\calS_{r+1}$ relative to $\tilde{\calS}_r$. From (\ref{et2}) we see that $e_t$ has non-vanishing parts only in the right cosets $\tilde{\calS}_r \cdot (k\,,\,r+1)$ and $\tilde{\calS}_r$. This leads to $h_s \not= e_t$ and $h_s \cdot e_t = e_t \not= h_s$.

Finally, we consider the case $k = r+1$. If $t$ is the tableau (\ref{tabl2}) with $k = r+1$, then $t' := (l , r+1) \circ t$ is a Young tableau\footnote{A Young tableau $t$ of $\calS_r$ can be regarded a one-to-one mapping of the boxes of the Young frame of $t$ onto the set $\{ 1 , \ldots , r \}$. $t$ maps every box onto that number which was placed into the box. Then the composition $p \circ t$ of a Young tableau $t$ and a permutation $p \in \calS_r$ is a Young tableau again.} (\ref{tabl2}) with $l = r+1$. A relation $t' = p \circ t$, $p \in \calS_{r+1}$, between Young tableaux leads to
\[
\calH_{t'} = p \circ \calH_t \circ p^{-1}
\;\;\;,\;\;\;
\calV_{t'} = p \circ \calV_t \circ p^{-1}
\;\;\;,\;\;\;
y_{t'} = p \cdot y_t \cdot p^{-1} \,.
\]
For the above tableaux $t' := (l , r+1) \circ t$ we obtain $y_{t'} = - (l , r+1) \cdot y_t$ since
$y_t \cdot (l , r+1)^{-1} = - y_t$.

Now, if we assume $h_s \cdot e_t = h_s$, then it follows from Lemma \ref{lem3}
\[
e_t \cdot h_s = e_t 
\;\;\Rightarrow\;\;
y_{t'} \cdot h_s = y_{t'}
\;\;\Rightarrow\;\;
e_{t'} \cdot h_s = e_{t'}
\;\;\Rightarrow\;\;
h_s \cdot e_{t'} = h_s \,.
\]
However, the last equation is a contradiction to our proof in the case $l = r+1$.
\end{proof}
\vspace{0.5cm}
\section{Use of $\nabla S$ and $\nabla A$ in generator formulas of algebraic covariant derivative curvature tensors}
Now we return to the question whether tensors (\ref{das}) can be used as generators $U$ for algebraic covariant derivative curvature tensors in formulas (\ref{divgen}). In \cite{fie22} we proved
\begin{Thm}
Let us denote by $S, A \in \calT_2 V$ symmetric or alternating tensors of order {\rm 2} and by $U \in \calT_3 V$ covariant tensors of order {\rm 3} whose symmetry class $\calT_{\frakr}$ is defined by a fixed minimal right ideal $\frakr$ from the equivalence class characterized by $(2 , 1) \vdash 3$. We consider the following types $\tau$ of tensors
\begin{eqnarray}
\tau: & &
\begin{array}{cccc}
y_{t'}^{\ast} (S \otimes U) & , & y_{t'}^{\ast} (U \otimes S) & , \\
y_{t'}^{\ast} (A \otimes U) & , & y_{t'}^{\ast} (U \otimes A) & , \\
\end{array} \label{types}
\end{eqnarray}
where $y_{t'} \in\bbK [\calS_5]$ is the Young symmetrizer of the standard tableau
\begin{eqnarray*}
t' & = &
\begin{array}{|c|c|c|c}
\cline{1-3}
1 & 3 & 5 & \\
\cline{1-3}
2 & 4 &\multicolumn{2}{c}{\;\;\;} \\
\cline{1-2}
\end{array}
\,.
\end{eqnarray*}
Then for each of the  above types $\tau$ the following assertions are equivalent:
\begin{enumerate}
\item{The vector space of algebraic covariant derivative curvature tensors $\fR' \in \calT_5 V$ is the set of all finite sums of tensors of the type $\tau$ considered. \label{statement1}}
\item{The right ideal $\frakr$ is different from the right ideal
$\frakr_0 := f \cdot \bbK [\calS_3]$ with generating idempotent
\begin{eqnarray}
f \;:=\; \left\{ \frac{1}{2}\,(\id - (1 \,3)) - \frac{1}{6}\,y \right\}
& \;\;\;,\;\;\; &
y \;:=\; \sum_{p \in \calS_3} \mathrm{sign}(p)\,p \,. \label{ipf2}
\end{eqnarray}
}
\end{enumerate} \label{thm3}
\end{Thm}
\noindent From Theorem \ref{thm3} we obtain easily
\begin{Thm}
Let $\nabla$ be a torsion-free covariant derivative on the mannifold $M$ and $p \in M$. Then Statement {\rm (\ref{statement1})} of Theorem {\rm \ref{thm3}} holds for the algebraic covariant derivative curvature tensors $\fR' \in \calT_5 M_p$ if we take the tensors $U$ from one of the tensor sets
\begin{eqnarray}
U \;=\; h_s^{\ast}(\nabla S)|_p & = & \nabla S|_p \,-\, \mathrm{sym}(\nabla S)|_p
\end{eqnarray}
or
\begin{eqnarray}
U \;=\; h_a^{\ast}(\nabla A)|_p & = & \nabla A|_p \,-\, \mathrm{alt}(\nabla A)|_p \;=\; \nabla A|_p \,-\, \mathrm{d}A|_p
\end{eqnarray}
formed from the whole of symmetric or alternating tensor fields $S , A \in \calT_2 M$.
\end{Thm}
\begin{proof}
Let us denote by $\calT_e$ the symmetry class $\calT_{\frakr}$ that is defined by a right ideal $\frakr = e \cdot \bbK [\calS_r]$ with generating idempotent $e \in \bbK [\calS_r]$. Then Proposition \ref{prop2} yields
\begin{eqnarray*}
\calT_{e_s} & = & \left\{ (\nabla S)|_p \;|\; S \in \calT_2 M \;\;\mathrm{symmetric} \right\} \\
\calT_{e_a} & = & \left\{ (\nabla A)|_p \;|\; A \in \calT_2 M \;\;\mathrm{alternating} \right\}
\end{eqnarray*}
from which we obtain\footnote{Note that the idempotents $e_s, f_s, h_s, e_a, f_a, h_a$ fulfil
$e_s^{\ast} = e_s$,
$f_s^{\ast} = f_s$,
$h_s^{\ast} = h_s$,
$e_a^{\ast} = e_a$,
$f_a^{\ast} = f_a$,
$h_a^{\ast} = h_a$.}
\begin{eqnarray*}
\calT_{h_s} & = &
\left\{ h_s (\nabla S)|_p \;|\; S \in \calT_2 M \;\;\mathrm{sym.} \right\}
\;=\;
\left\{ \nabla S|_p \,-\, \mathrm{sym}(\nabla S)|_p \;|\; S \in \calT_2 M \;\;\mathrm{sym.} \right\} \\
\calT_{h_a} & = & \left\{ h_a (\nabla A)|_p \;|\; A \in \calT_2 M \;\;\mathrm{alt.} \right\}
\;=\;
\left\{ \nabla A|_p \,-\, \mathrm{alt}(\nabla A)|_p \;|\; A \in \calT_2 M \;\;\mathrm{alt.} \right\}\,.
\end{eqnarray*}

Now we have only to check that the idempotents
\begin{eqnarray*}
h_s \;=\; e_s \,-\, f_s
& \;\;\;,\;\;\; &
h_a \;=\; e_a \,-\, f_a
\end{eqnarray*}
do not generate the right ideal $\frakr_0$. We do this by verifying
\begin{eqnarray}
f \cdot h_s \;\not=\; h_s
& \;\;\;,\;\;\; &
f \cdot h_a \;\not=\; h_a \,. \label{check}
\end{eqnarray}
The fastes way would be a computer calculation by means of {\sf PERMS} \cite{fie10}. A calculation by hand has the starting point
(\ref{eee}), (\ref{idpotf}), (\ref{ipf2}) and $z := \frac{1}{2}(\id - (1\,3))$.
From the rules
\begin{itemize}
\item{''symmetrization + alternation = 0''}
\item{''alternation + alternation = alternation''}
\end{itemize}
we obtain immediately
\begin{eqnarray*}
y \cdot f_s \;=\; 0
& \;\;\;,\;\;\; &
y \cdot e_s \;=\; 0
\;\;\;\;,\;\;\;\; 
z \cdot f_s \;=\; 0 \\
y \cdot f_a \;=\; y
& \;\;\;,\;\;\; &
y \cdot e_a \;=\; y
\;\;\;\;,\;\;\;\; 
z \cdot f_a \;=\; f_a
\end{eqnarray*}
Furthermore we have the products
\begin{eqnarray*}
z \cdot e_s  & = &
\frac{1}{4} \left\{ [1,2,3] + [2,1,3] - [2,3,1] - [3,2,1] \right\} \\
z \cdot e_a & = &
\frac{1}{4} \left\{ [1,2,3] - [2,1,3] + [2,3,1] - [3,2,1] \right\} \,.
\end{eqnarray*}
This leads to
\begin{eqnarray*}
f \cdot h_s & = & z \cdot e_s \;=\;
\frac{1}{4} \left\{ [1,2,3] + [2,1,3] - [2,3,1] - [3,2,1] \right\}
\end{eqnarray*}
and
\begin{eqnarray*}
f \cdot h_a & = &
z \cdot e_a - z \cdot f_a - \frac{1}{6}\, y \cdot e_a + \frac{1}{6}\, y \cdot f_a \\
& = &
z \cdot e_a - f_a - \frac{1}{6}\, y  + \frac{1}{6}\, y  \\
& = &
z \cdot e_a - f_a \,.
\end{eqnarray*}
But we see from these results that (\ref{check}) is valid because
\begin{itemize}
\item{$h_s$ is a linear combination of 6 permutations and $z \cdot e_s$ has only 4 summands,}
\end{itemize}
\hspace*{28pt}$\bullet$ $z \cdot e_a \not= e_a$.
\end{proof}
\vspace{1cm}
\section*{Appendix}
Using the {\sf Mathematica} package {\sf PERMS} \cite{fie10} for $r = 2, 3, 4$ we found Young symmetrizer idempotents $e_t = {\mu}_t y_t$ with Young frames $(2 , 1^{r-1}) \vdash r + 1$ which reproduce ($h_a \cdot e_t = h_a$) or annihilate ($h_a \cdot e_t = 0$) the idempotent $h_a \in \bbK [\calS_{r+1}]$ of the symmetric group $\calS_{r+1}$ considered. Here we present complete lists of the Young tableaux $t$ for which the $e_t$ possess such a property.

\subsubsection*{$2$ tableaux for $r = 2$ such that $h_a \cdot e_t = h_a$}
\begin{verbatim}
  {1, 3}, {2, 3}
  {2}     {1}
\end{verbatim}

\subsubsection*{$2$ tableaux for $r = 2$ such that $h_a \cdot e_t = 0$}
\begin{verbatim}
  {1, 2}, {2, 1}
  {3}     {3}
\end{verbatim}

\subsubsection*{$6$ tableaux for $r = 3$ such that $h_a \cdot e_t = h_a$}
\begin{verbatim}
  {1, 4}, {1, 4}, {2, 4}, {2, 4}, {3, 4}, {3, 4}
  {2}     {3}     {1}     {3}     {1}     {2}
  {3}     {2}     {3}     {1}     {2}     {1}
\end{verbatim}

\subsubsection*{$12$ tableaux for $r = 3$ such that $h_a \cdot e_t = 0$}
\begin{verbatim}
  {1, 2}, {1, 2}, {1, 3}, {1, 3}, {2, 1}, {2, 1}, {2, 3}, {2, 3}, 
  {3}     {4}     {2}     {4}     {3}     {4}     {1}     {4}
  {4}     {3}     {4}     {2}     {4}     {3}     {4}     {1}
 
  {3, 1}, {3, 1}, {3, 2}, {3, 2}
  {2}     {4}     {1}     {4}
  {4}     {2}     {4}     {1}
\end{verbatim}

\subsubsection*{$24$ tableaux for $r = 4$ such that $h_a \cdot e_t = h_a$}
\begin{verbatim}
  {1, 5}, {1, 5}, {1, 5}, {1, 5}, {1, 5}, {1, 5}, {2, 5}, {2, 5}, 
  {2}     {2}     {3}     {3}     {4}     {4}     {1}     {1}
  {3}     {4}     {2}     {4}     {2}     {3}     {3}     {4}
  {4}     {3}     {4}     {2}     {3}     {2}     {4}     {3}
 
  {2, 5}, {2, 5}, {2, 5}, {2, 5}, {3, 5}, {3, 5}, {3, 5}, {3, 5}, 
  {3}     {3}     {4}     {4}     {1}     {1}     {2}     {2}
  {1}     {4}     {1}     {3}     {2}     {4}     {1}     {4}
  {4}     {1}     {3}     {1}     {4}     {2}     {4}     {1}
 
  {3, 5}, {3, 5}, {4, 5}, {4, 5}, {4, 5}, {4, 5}, {4, 5}, {4, 5}
  {4}     {4}     {1}     {1}     {2}     {2}     {3}     {3}
  {1}     {2}     {2}     {3}     {1}     {3}     {1}     {2}
  {2}     {1}     {3}     {2}     {3}     {1}     {2}     {1}
\end{verbatim}

\subsubsection*{$72$ tableaux for $r = 4$ such that $h_a \cdot e_t = 0$}
\begin{verbatim}
  {1, 2}, {1, 2}, {1, 2}, {1, 2}, {1, 2}, {1, 2}, {1, 3}, {1, 3}, 
  {3}     {3}     {4}     {4}     {5}     {5}     {2}     {2}
  {4}     {5}     {3}     {5}     {3}     {4}     {4}     {5}
  {5}     {4}     {5}     {3}     {4}     {3}     {5}     {4}
 
  {1, 3}, {1, 3}, {1, 3}, {1, 3}, {1, 4}, {1, 4}, {1, 4}, {1, 4}, 
  {4}     {4}     {5}     {5}     {2}     {2}     {3}     {3}
  {2}     {5}     {2}     {4}     {3}     {5}     {2}     {5}
  {5}     {2}     {4}     {2}     {5}     {3}     {5}     {2}
 
  {1, 4}, {1, 4}, {2, 1}, {2, 1}, {2, 1}, {2, 1}, {2, 1}, {2, 1}, 
  {5}     {5}     {3}     {3}     {4}     {4}     {5}     {5}
  {2}     {3}     {4}     {5}     {3}     {5}     {3}     {4}
  {3}     {2}     {5}     {4}     {5}     {3}     {4}     {3}
 
  {2, 3}, {2, 3}, {2, 3}, {2, 3}, {2, 3}, {2, 3}, {2, 4}, {2, 4}, 
  {1}     {1}     {4}     {4}     {5}     {5}     {1}     {1}
  {4}     {5}     {1}     {5}     {1}     {4}     {3}     {5}
  {5}     {4}     {5}     {1}     {4}     {1}     {5}     {3}
 
  {2, 4}, {2, 4}, {2, 4}, {2, 4}, {3, 1}, {3, 1}, {3, 1}, {3, 1}, 
  {3}     {3}     {5}     {5}     {2}     {2}     {4}     {4}
  {1}     {5}     {1}     {3}     {4}     {5}     {2}     {5}
  {5}     {1}     {3}     {1}     {5}     {4}     {5}     {2}
 
  {3, 1}, {3, 1}, {3, 2}, {3, 2}, {3, 2}, {3, 2}, {3, 2}, {3, 2}, 
  {5}     {5}     {1}     {1}     {4}     {4}     {5}     {5}
  {2}     {4}     {4}     {5}     {1}     {5}     {1}     {4}
  {4}     {2}     {5}     {4}     {5}     {1}     {4}     {1}
 
  {3, 4}, {3, 4}, {3, 4}, {3, 4}, {3, 4}, {3, 4}, {4, 1}, {4, 1}, 
  {1}     {1}     {2}     {2}     {5}     {5}     {2}     {2}
  {2}     {5}     {1}     {5}     {1}     {2}     {3}     {5}
  {5}     {2}     {5}     {1}     {2}     {1}     {5}     {3}
 
  {4, 1}, {4, 1}, {4, 1}, {4, 1}, {4, 2}, {4, 2}, {4, 2}, {4, 2}, 
  {3}     {3}     {5}     {5}     {1}     {1}     {3}     {3}
  {2}     {5}     {2}     {3}     {3}     {5}     {1}     {5}
  {5}     {2}     {3}     {2}     {5}     {3}     {5}     {1}
\end{verbatim}
\newpage
\begin{verbatim}
  {4, 2}, {4, 2}, {4, 3}, {4, 3}, {4, 3}, {4, 3}, {4, 3}, {4, 3}
  {5}     {5}     {1}     {1}     {2}     {2}     {5}     {5}
  {1}     {3}     {2}     {5}     {1}     {5}     {1}     {2}
  {3}     {1}     {5}     {2}     {5}     {1}     {2}     {1}
\end{verbatim}
\vspace{0.5cm}
\noindent {\bf Acknowledgement.} I would like to thank P. B. Gilkey for suggesting these investigations and for many important and helpful discussions.
\vspace{0.4cm}

\end{document}